\documentclass[fleqn,10pt]{article}
\usepackage{graphicx}
\usepackage{latexsym}
\usepackage{amsmath}
\usepackage{amssymb}
\usepackage{amsthm}
\usepackage{hyperref}
\def \Q {{\mathbb Q}}

\def \Z {{\mathbb Z}}

\def \C {{\mathbb C}}

\def \al {{\alpha}}

\def \ze {{\zeta}}
\newtheorem{theorem}{Theorem}[section]

\newtheorem{lemma}[theorem]{Lemma}

\newtheorem{rem}[theorem]{Remark}

\author{A. Satyanarayana Reddy\footnote{Department of Mathematics and
Statistics, Indian
Institute of Technology Kanpur,  India 208016; (e-mail:
satya@iitk.ac.in).}}

\title{The lowest degree $0,1$-polynomial divisible by cyclotomic polynomial}
\date{}
\begin{document}
\maketitle

%%%%%%%%%%%%%%%%%%%%%%%%%%%%%%%%%%%%%%%%%%%%%%%%%%%%%%%%%%%%%%%%%%%
\begin{abstract}
Let $n$ be an even positive integer with at most three distinct
prime factors and let $\ze_n$ be a primitive $n$-th root of unity.
In this study, we made an attempt to find  the lowest-degree
$0,1$-polynomial $f(x) $ having at least three terms such
that $f(\ze_n)$ is a minimal vanishing sum of distinct $n$-th roots of
unity.
\end{abstract}
%%%%%%%%%%%%%%%%%%%%%%%%%%%%%%%%%%%%%%%%%%%%%%%%%%%%%%%%%%%%%%

\section{Introduction and preliminaries}

For a fixed positive integer $n$, let $U_n = \{ k : 1 \le k \le n,
\gcd(k,n) = 1\}$. If $|S|$ denotes the cardinality of the set
$S$, then $|U_n| = \varphi(n)$, the well known {\em Euler-totient
function}.  Let $\zeta_n\in\C$ denote a primitive $n$-th root of
unity, {\it i.e.,} $\zeta^n=1$ and $\zeta^k\neq 1$ for $0<k<n$. It
is well known that the  {\em $n$-th cyclotomic polynomial},
$\Phi_n(x)=\prod\limits_{k \in U_n}(x-\zeta_n^k)\in\Z[x],$
 is the minimal polynomial of $\ze_n$ and  $\deg(\Phi_n(x)) = \varphi(n)$.
 It is also known that
$x^n-1=\prod\limits_{d \mid n}\Phi_d(x)$ (here $a \mid b$ means `$a$ divides
$b$'). Now, using the property of minimal polynomials it follows
that  whenever $f(\ze_n) = 0$ for some $f(x) \in \Q[x]$ then $\Phi_n(x)$ divides
$f(x)$.  
 Let $\Phi_n(x)=\sum\limits_{k=0}^{\varphi(n)}a_{n,k} x^k$. 
Then it can be easily verified  that $a_{n,i}=a_{n,\varphi(n)-i}$, where $0\leq
i\leq \varphi(n)$,  $a_{n,0}=1$ (for $n>1$) and $a_{n,1}=-\sigma(n)$, where
$\sigma(n) = \sum\limits_{k \in U_n} \zeta_n^k$, is the sum of primitive $n$-th
roots of unity.
Before stating a result that gives the value of $a_{n,1}$, recall that a
positive
integer is said to be {\em square-free} if its decomposition into prime
numbers/factors does not have any repeated factors.

\begin{lemma}[Apostol~\cite{T:A}]\label{lem:mu}
Let $n$ be a positive integer. Then $\sigma(n)= \mu(n)$, where
$$ \mu(n)=
\begin{cases} 0, & {\mbox{ if }} \; n \; {\mbox{ is not square-free,}} \\
 1, &  {\mbox{ if }} \; n \; {\mbox{ has even number of prime factors,}}\\
-1, & {\mbox{ if }} \; n \; {\mbox{ has odd number of prime factors.}}
\end{cases}$$
\end{lemma}

For  $f(x)\in \Z[x]$,  we denote  the set of
 coefficients of $f(x)$ by $V_f$
  and
the set of exponents of $f(x)$ by $E_f$ or $E_{f(x)}$ .  We define the {\em height} of the
polynomial
$f(x)\in\Z[x]$ to be the maximum absolute value of a coefficient
of $f(x)$. Let $A(n)$ be the height of $\Phi_n(x)$. We say that a
cyclotomic polynomial $\Phi_n(x)$ is {\em flat} if $A(n) = 1$.
For $n < 105$, $\Phi_n(x)$ is flat. It was once conjectured that
all cyclotomic polynomials are flat, however $A(105) = 2$. In
fact, $A(n)$ is unbounded refer Emma Lehmer~\cite{lehmer}.  It is known that
$\Phi_{pn}(x)=\Phi_n(x^p)$ whenever $p|n$ and $\Phi_{2n}(x)=\Phi_n(-x)$ whenever $n$ is odd.
Hence, if  $\Phi_n(x)$ is a flat polynomial, then so is  $\Phi_{pn}(x)$, where $p|n$ and  $\Phi_{2n}(x)$ where $n$ is odd.

Let $k$ be the number of distinct odd
prime factors of $n$. For square-free $n$, this number $k$ is
called the  {\em order} of the cyclotomic polynomial
$\Phi_n(x)$.  It is known that all cyclotomic polynomials of order
$1$ and order $2$ are flat.
Gennady Bachman~\cite{gb} gave the first infinite family of flat cyclotomic
polynomials of order three
 and this family was expanded by Kaplan~\cite{kap1}. In \cite{kap2} Kaplan gave
some flat
polynomials of order four. It is unknown whether there
are any flat cyclotomic polynomials of order greater than four.

Now, recall that an equation  of the form $\ze_n^{k_1}+\ze_n^{k_2}+\dots
+\ze_n^{k_l}=0$ is called a
 {\em vanishing sum} of  $n$-th roots of unity of weight $l$, 
 where $k_1,k_2,\ldots,k_l\in \{0,1,2,\ldots, n-1\}$.  The vanishing sum is said
to be {\em minimal}
if no proper sub-sums thereof can be zero. 

For example, Let  $p$ be a prime divisor of $n$.
Then $\zeta_n^{(n/p) k}$ for $1 \le k \le p-1$ are primitive $p$-th roots of
unity and hence using Lemma~\ref{lem:mu}, it follows that $1+\sigma(p)=0$ is a
vanishing sum of $n$-th roots of unity and is in fact minimal as well.
\begin{rem}
A vanishing sum can always be multiplied by a root of unity to get
another vanishing sum. As such, the two vanishing sums are similar and in
literature, it is said that the latter is obtained from the former by a
{\em rotation}. Hence, the classification of minimal vanishing sums needs to be done
only up to
rotations (by roots of unity).
\end{rem}

We now state Corollary~$3.2$ of T.Y.Lam and K.H.Leung~\cite{Lam1} that relates
vanishing sums with square-free positive integers.
\begin{lemma}\label{lem:n0}
 If $\al_1+\al_2+\dots +\al_m=0$ is a  minimal vanishing sum of $n$-th roots of
unity, then after a suitable rotation, we may assume that all
$\al_i$'s are $n_0$-th roots of unity where $n_0$ is the largest
square-free part of $n$.
\end{lemma}

For example, if $n=p^k$ for some prime $p$, then $1+\sigma(p)=0$ is the  only
minimal vanishing sum of roots of unity  up to rotation.

\begin{lemma}\label{lem:two primes}[T.Y.Lam and
K.H.Leung~\cite{Lam1}]
Let $n=p^a q^b$, where $p,q$ are distinct primes. Then, up to a
rotation, the only minimal vanishing sums of $n$-th roots of unity
are: $1+\sigma(p)=0$ and $1+\sigma(q)=0$.
\end{lemma}

 If $n$ has 3 or more distinct  prime factors, then 
T.Y.Lam and K.H.Leung~\cite{Lam1} provided another minimal
vanishing sum of $n$-th roots of unity, $\sigma(p_i)\cdot
\sigma(p_j)+\sigma(p_k)=0$ where $p_i,p_j,p_k$ are distinct primes
dividing $n$. Finding all minimal vanishing sums of $n$-th roots
of unity when $n$ has 3 or more distinct prime divisors seems difficult.
In Section~\ref{sec:two},  we constructed another minimal vanishing sum  
of $n$-th roots of unity, when $n$ is even and $\Phi_n(x)$ is flat. And we use
this minimal vanishing sum for finding lowest $0,1$ polynomial  divisible by
$\Phi_n(x)$.
Note that  in their
paper, T.Y.Lam and K.H.Leung~\cite{Lam1}, allowed the repetition
of $n$-th roots of unity in the vanishing sum. Gary
Sivek~\cite{sivek} considered the vanishing sums of distinct
$n$-th roots of unity.
He found the number of terms in those vanishing sums. Unless
specified otherwise, from now onwards, 
vanishing sum means vanishing sum of distinct
$n$-th roots of unity.

John Steinberger~\cite{jsb} studied the lowest-degree polynomial
with nonnegative coefficients  divisible by $\Phi_n(x)$. He
conjectured that the lowest degree monic polynomial with nonnegative
coefficients divisible by $\Phi_n(x)$, $n > 1$, is $1 +
x^{n/p}+\dots+ x^{(p-1)n/p}$ where $p$ is the smallest prime
dividing $n$. And proved this conjecture  holds when $n$ is even
or when $n$ is a prime power or when $2/p >
1/q_1+\dots + 1/q_k$ where $q_1,\ldots, q_k$ are the other primes
besides $p$ dividing $n$. Consequently, if $n$ is a  odd prime power
or when $2/p > 1/q_1+\dots + 1/q_k$ where $q_1,\ldots, q_k$ are
the other primes besides $p$ dividing $n$ and p is odd, then  the lowest-degree  $0,1$-
polynomial divisible by $\Phi_n(x)$ having at least three terms is $1 +
x^{n/p}+\dots+ x^{(p-1)n/p}$.

Let $n=2m$ for some $m\in \Z$. Then  the lowest-degree
$0,1$-polynomial divisible by $\Phi_n(x)$ is $x^{n/2}+1$.  
Our interest in this paper, is to  find  the
lowest-degree $0,1$-polynomial $f(x)$ with at least three terms such
that $f(\ze_{2m})$ is a minimal vanishing sum of $2m$-th roots of
unity.  It is  known  (for example, see John P. Steinberger~\cite{J:S}) that
the problem of  finding the polynomials divisible by $\Phi_n(x)$ is
equivalent to finding polynomials divisible by $\Phi_{n_0}(x),$
where $n_0$ is the maximum square-free factor of $n$. With this observation and Lemma~\ref{lem:n0},
it is sufficient to consider square free $n$. Hence, in the remaining part of this paper, we consider
$n=2p_1\dots p_k$, where $k\geq 1$  and $p_1<p_2\dots <p_k$ are prime numbers.

Let us denote $$ G_n=\{f(x)| f(\ze_n)\mbox{ is a minimal vanishing
sum}\}.$$  Consequently, by our assumption of vanishing sum, if $f(x)\in G_n$, then $V_f\subseteq \{0,1\}$, $\deg(f(x))\leq n-1$ and $\Phi_n(x)|f(x)$.
Recall that  a vanishing sum $f_1(\ze_n)=0$ is a rotation
of another vanishing sum $f_2(\ze_n)=0$ if $f_1(\ze_n)=\ze_n^t
f_2(\ze_n)$ for some $t$. Now we define a relation $\mathcal{R}$ on
$G_n$ as $f_1(x)\mathcal{R} f_2(x)$ if and only if $f_1(\ze_n)=0$ is
a rotation of $f_2(\ze_n)=0$.
 Clearly $\mathcal{R}$ is an equivalence relation on $G_n$. Let
$H_n=G_n/\mathcal{R}$ be the set of equivalence classes of $G_n$
under the equivalence relation $\mathcal{R}$. Without loss of
generality we assume that the elements of $H_n$ are the polynomials
of least degree in the corresponding equivalence classes. Under this assumption
$H_n\subseteq G_n$. Hence, if $f(x)\in H_n$, then $V_f\subseteq \{0,1\}$,
$f(\ze_n)$  is a minimal vanishing sum, $\deg(f(x))\leq n-1$ and $f(0)=1$.
For example, let $n=p_0p_1\dots p_k,\;p_0<p_1<\dots <p_k\;
,\;p_0=2$. Let
$g_i(x)=1+\sum_{k=1}^{p_i-1}x^{\frac{kn}{p_i}}=1+x^{\frac{n}{p_i}}+x^{\frac{2n}{
p_i}}\dots
+x^{\frac{(p_i-1)n}{p_i}}$.  Then $g_i(x)$ is the polynomial
corresponding to minimal vanishing sum $1+\sigma(p_i)=0$. Hence $g_i(x)\in G_n$. It is also easy to see that  $g_i(x)\in H_n$ and $g_0(x)=x^{n/2}+1$ is the only polynomial in $H_n$ with two terms.
With the definition of $H_n$, the problem of our interest is to find {\it a} lowest degree polynomial in $H_n$ with at least three terms.

%%%%%%%%%%%%%%%%%%%%%%%%%%%%%%%%%%%%%%%%%%%%%%%%%%%%%%%%%%%%%%%%%%%
\section{{\small{Minimal vanishing sum from flat cyclotomic
polynomial}}}\label{sec:two}
%%%%%%%%%%%%%%%%%%%%%%%%%%%%%%%%%%%%%%%%%%%%%%%%%%%%%%%%%%%%%%%%%%%%

 Let $n=p_0p_1p_2\dots p_k$, where $k\geq
1$ and $p_0=2<p_1<p_2<\dots < p_k$ are distinct primes.
Suppose $\Phi_n(x)$ is flat.  One can verify that
\begin{equation}\label{eq:Phi}
\Phi_n(x)= \begin{cases} x^{\varphi(n)}-x^{\varphi(n)-1}\pm \dots -x+1, &
{\mbox{ when }} \; k \; {\mbox{ odd,}} \\
x^{\varphi(n)}+x^{\varphi(n)-1}-x^{\varphi(n)-p_1}\pm \dots
-x^{p_1}+x+1, & {\mbox{ when }} \; k \; {\mbox{  even. }}
\end{cases}
\end{equation}
 We can write $\Phi_n(x)=f_1(x)-f_2(x)$, where
$f_1(x)$ and $f_2(x)$ are polynomials with positive coefficients.
Observe this representation is unique as  $E_{f_1}\cap E_{f_2}=\emptyset$.
Now
\begin{eqnarray*}
\Phi(\ze_n)&=&0\\
&\Rightarrow & f_1(\ze_n)-f_2(\ze_n)=0\\
&\Rightarrow & f_1(\ze_n)+\ze_n^{n/2}\cdot f_2(\ze_n)=0\\
&\Rightarrow & \Phi_n(x)| (f_1(x)+x^{n/2}f_2(x)).
\end{eqnarray*}

Let $\Phi^T_n(x)= f_1(x)+x^{n/2}f_2(x)$. Then $\Phi^T_n(x)$ is a
$0,1$-polynomial.
And from Equation~(\ref{eq:Phi}) we have 
\begin{equation}
 \deg(\Phi^T_n(x))=\begin{cases} \phi(n)-1+\frac{n}{2},\;\; \mbox{when\; k\;
is\;
odd},\\
\phi(n)-p_1+\frac{n}{2},\;\; \mbox{when \; k\; is\; even.}
                    \end{cases}
\end{equation}

\begin{lemma}\label{lem:flat}
 Let $k\geq 1$ and  $n=p_0p_1p_2\dots p_k$, where $p_0=2<p_1<p_2<\dots < p_k$
 are distinct primes.  Suppose  $\Phi_n(x)$ is
a flat polynomial. Then $\Phi_n^T(\ze_n)=0$ is a minimal vanishing
sum of $n$-th roots of unity.
\end{lemma}

\begin{proof}
 Let us suppose $\Phi_n^T(\ze_n)=0$ is not a minimal
vanishing sum of $n$-th roots of unity. Suppose $A(\ze_n)$ is a
proper sub-sum of  $\Phi_n^T(\ze_n)$ such that $A(\ze_n)=0$. Now
we write $\Phi_n^T(x)=A(x)+B(x)$. Consequently $B(\ze_n)=0$.
Without loss of generality assume that  $A(x)$
 does not contain the term $x^{\varphi(n)}$.
We now  write $A(x)=A_1(x)+A_2(x)$ such that the exponent of every
term in $A_1(x)>n/2$ and that of $A_2(x)<n/2$. Then
$\Phi_n(x)|A_2(x)-x^{-n/2}A_1(x)$ which is not possible since
$\deg(A_2(x)-x^{-n/2}A_1(x))< \varphi(n)=\deg(\Phi_n(x))$.
\end{proof}

The following result is a direct consequence of Lemma~\ref{lem:two primes}. We
are giving the proof for the sake of completeness.

\begin{lemma}
Let $p$ be any positive odd prime integer and  $n=2p$, then
$H_n=\{g_0,g_1\}$. Further $g_1(x)=\Phi^T_n(x)$.
\end{lemma}

\begin{proof}
First observe that
$\Phi_{2p}(x)=\Phi_p(-x)=x^{p-1}-x^{p-2}+x^{p-3}-\dots-x+1$. Hence
$\Phi^T_n(x)=\sum_{i=0}^{p-1}x^{2i}=g_1(x)$. Let us suppose that
$f(x)=x^{k_1}+x^{k_2}+\dots +x^{k_l}+1\in H_n$ where
$n>k_1>k_2>\dots >k_l>0$ and $f(x)\neq g_0(x)$. By
definition of $H_n$,  $\frac{n}{2}\notin E_f$. Now  we write
$f(x)=f_1(x)+f_2(x)$ in such a way that $min(E_{f_1})>\frac{n}{2}$
and $max(E_{f_2})<\frac{n}{2}$. Let
$h(x)=f_2(x)-x^{\frac{n}{2}}f_1(x)$. Then degree of
$\deg(h(x))\leq \frac{n}{2}-1=p-1=\varphi(n)$. But
$h(\ze_n)=f_2(\ze_n)-\ze_n^{\frac{n}{2}}f_1(\ze_n)=f(\ze_n)=0$. Hence
$\Phi_n(x)|h(x)$ and $h(x)=\Phi_n(x)$. Thus $f(x)=\Phi_n^T(x)=g_1(x)$.
\end{proof}

From the above lemma, the lowest degree polynomial in $H_{2p}$ with at least three terms is $\Phi_{2p}^T(x)$. 
We now  try to find a lowest degree polynomial in $H_n$ with at least three terms, whenever
$n=2p_1p_2$, where  $p_1<p_2$ are odd primes. Since $g_1(x)\in H_n$,  the
degree of required polynomial is $\leq n-\frac{n}{p_1}$. In
1883, Miggotti~\cite{mig}, showed that $\Phi_n(x)$ is a flat polynomial.
Hence from Lemma~\ref{lem:flat}, $\Phi_n^T(x)\in G_n$. In the following
subsection, we show that $\Phi_n^T(x)\in H_n$.

%%%%%%%%%%%%%%%%%%%%%%%%
\subsection{$n=2pq$}
%%%%%%%%%%%%%%%%%%%%%%%
In this subsection, we suppose $n=2pq$ where
$p<q$ are odd primes. Lam and Heung~\cite{Lam2}
gave  a  nice expression for 
$\Phi_{pq}(x)$ as
\begin{equation}\label{eq:cyc}
 \Phi_{pq}(x)=(\sum_{i=0}^rx^{ip})(\sum_{j=0}^sx^{jq})-(\sum_{i=r+1}^{q-1}x^{ip}
)(\sum_{j=s+1}^{p-1}x^{jq})x^{-pq},
\end{equation}

where $r$ and $s$ are positive  integers such that
 $rp+sq=(p-1)(q-1)=\varphi(n),\;0\leq r\leq q-2\;and\;0\leq j\leq p-2$.
Hence $r$ and $s$  have same parity.

Let us write Equation~\ref{eq:cyc} as $\Phi_{pq}(x)=[A(x)+B(x)]-[C(x)+D(x)]$
such that  $A(x)+B(x)= (\sum_{i=0}^rx^{ip})(\sum_{j=0}^sx^{jq})$ and
 $C(x)+D(x)=
(\sum_{i=r+1}^{q-1}x^{ip})(\sum_{j=s+1}^{p-1}x^{jq})x^{-pq}$,
where \begin{eqnarray*}E_{A} &=& \{ip+jq| 0\leq i\leq r, 0\leq j\leq s\;\mbox{ and both $i$ and $j$
have same parity}\},\\
E_{B}&=&\{ip+jq| 0\leq i\leq r,0\leq j\leq s\;\mbox{ and $i$ and $j$ are of
different parity}\}.
\end{eqnarray*}
It is easy to see that $max(E_A)=\varphi(n)$, $min(E_A)=0$, $max(E_B)=\varphi(n)-p$ and $min(E_B)=p$. Similarly we can define $C(x)$ and $D(x)$ such that 
$A(x),\; C(x)$ are even functions and $B(x),\;D(x)$ are odd functions.

With this we can write
\begin{equation*}
 \Phi_{2pq}(x)=[A(x)+D(x)]-[B(x)+C(x)]\;\mbox{and}\;
\Phi_{2pq}^T(x)=[x^{n/2}(B(x)+C(x))]+[A(x)+D(x)].
\end{equation*}
In order to prove $\Phi^T_n(x)\in H_n$ it is sufficient to prove that
$\deg(\Phi^T_n(x))\leq \deg(x^i\Phi^T_n(x)) \;for\;all\; i\in
E_{\Phi^T_n}$. Hence the following result.

\begin{theorem}\label{thm:2pq}
 Let $n=2pq$ where $p,q (p<q)$ are odd primes. Then
\begin{enumerate}
 \item \label{thm:2pq:1} $\deg(x^i\Phi_n^T(x))>\deg(\Phi^T_n(x))\; for\;all\;\; i\ne 0,\frac{n}{2}-p$.
\item \label{thm:2pq:2}
$\deg(\Phi_n^T(x))=\deg(x^{\frac{n}{2}-p}\Phi_n^T(x))=\frac{n}{2}
+\varphi(n)-p=n-2p-q+1$.
\end{enumerate}
\end{theorem}
First we prove the following lemma.

\begin{lemma}
 Let $n=2pq$, where $p,q (p<q)$ are odd primes.
 Let $h\in E_{\Phi_{n}}$, $h<\varphi(n)$, then  there exists $g\in E_{\Phi_n}$ such
that $0<g-h\leq p+q$.
\end{lemma}

\begin{proof}
If we prove the result for $E_{A}$ and $E_{B}$, then it is
also true of $E_{A(x)+D(x)}$ and $E_{B(x)+C(x)}$. Hence the result
follows for $E_{\Phi_{n}}$  as $\deg(A(x))-\deg(B(x))=p$.

 First we write $A(x)=A_1(x)+A_2(x)$ such that
 \begin{eqnarray*}
E_{A_1} &=& \{ip+jq| 0\leq i\leq r,
0\leq j\leq s\;and\;both\;i\;and\; j \;are \;even\},\\
 E_{A_2} &=& \{ip+jq| 0\leq i\leq r,
0\leq j\leq s\;and\;both\;i\;and\; j \;are \;odd\} .
\end{eqnarray*}

Let $h\in E_{A}$. We will find  $g\in E_{A}$ such that $g>h$ and
$g-h\leq p+q$.

First consider the case when  both $r$ and $s$ are even.\\
Suppose $h\in E_{A_1}$. Then $h=ip+jq$ where $i,j$ are even,
$0\leq i\leq r$ and $0\leq j\leq s$. Also $i=r, j=s$ is not possible since
$h=rp+sq=\varphi(n)$.

\begin{itemize}
\item[Case1:] $i\leq r-2,j\leq s-2 \Rightarrow g=(i+1)p+(j+1)q\in E_{A_2}$.
\item[Case2:] $i=r\;and\; j\leq s-2 \Rightarrow g=(r-1)p+(j+1)q\in E_{A_2}$.
\item[Case3:] $i\leq r-2,j=s \Rightarrow g=(i+2)p+sq\in E_{A_1}$.
\end{itemize}

Now suppose $h\in E_{A_2}$. Then $h=ip+jq$, $i,j$ are odd,
$0\leq i\leq r-1$ and $0\leq j\leq s-1$. In this case we can choose
$g=(i+1)p+(j+1)q\in E_{A_1}$.  With the same arguments we can prove that the
difference between any two
consecutive elements in $E_{A}(\mbox{with an abuse of language, assume  there is an ordering in $E_f$} )$ is at most $p+q$ when both $r$
and $s$ are odd.

Now we prove that the difference between any two consecutive elements in
 $E_{B}$ is at most $p+q$. Before proceeding, recall $\varphi(n)-p$ is the
largest element in $E_{B}$. Hence we suppose that $h<\varphi(n)-p$. As if
$h=\varphi(n)-p$, then we choose $g=\varphi(n)\in E_A$ and $g-h<p+q$.

Let $h\in E_B$. Then $h=ip+jq$, where $i,j$ have different parity and 
$0\leq i\leq r$ and $0\leq j\leq s$. Also by our assumption of $h$, $i=r-1,j=s$ is not possible.
We have the following cases.
\begin{itemize}
 \item [Case1:] $i$ and $r$ have same parity. Then $j$ and $s$ have different parity.
\begin{enumerate}
\item $i\leq r-2,\;j\leq s-1\Rightarrow g=(i+1)p+(j+1)q$.
\item $i=r\;and\; j\leq s-1 \Rightarrow g=(r-1)p+(j+1)q$.
\end{enumerate}
 \item [Case2:] $i$ and $r$ have different parity. Then $j$ and $s$ have same parity.
\begin{enumerate}
\item $i\leq r-1,\;j\leq s-2\Rightarrow g=(i+1)p+(j+1)q$.
\item $i<r-1\;and\; j=s \Rightarrow g=(i+2)p+sq$.
\end{enumerate}
\end{itemize}

\end{proof}

We now prove the Theorem~\ref{thm:2pq}.

\begin{proof} Proof of Part~\ref{thm:2pq:1}:
From the above lemma, the difference  between any 
two consecutive terms in $E_{x^{\frac{n}{2}}(C(x)+B(x))}$ and $E_{A(x)+D(x)}$ 
is at most $p+q$. It is  easy to see that the terms with the largest and the smallest exponent in
$A(x)+D(x)$ belong to $A(x)$. Similar statement hold for $B(x)$ in
$B(x)+C(x)$. Hence we have 
$\underbrace{min(E_{x^{\frac{n}{2}}(C(x)+B(x))})-max(E_{A(x)+D(x)})}=
\frac{n}{2}+p-\varphi(n)=\underbrace{n-max(E_{x^{\frac{n}{2}}(C(x)+B(x))})}$.
With these observations, the largest difference in
the exponents of two consecutive terms in $\Phi_n^T(x)$ is
$\frac{n}{2}+p-\varphi(n)=2p+q-1(>p+q)$.
  Thus we conclude that
$\deg(x^i\Phi_n^T(x))>\deg(\Phi_n^T(x))$  for all  $i\ne 0,\frac{n}{2}-p$. \\
Proof of Part~\ref{thm:2pq:2}: It is easy to see.
\end{proof}

Now the degree of $\Phi_n^T(x)$ is
$\frac{n}{2}+\varphi(n)-p=n-2p-q+1$, whereas the degree of
$g_1(x)=n-\frac{n}{p}=n-2q$, where $g_1(x)$ is the polynomial corresponding to minimal vanishing sum $1+\sigma(p)=0$. 
Now we conclude this work with the following conjecture and we
feel that Lemma~\ref{lem:s} given below, may be useful in proving this conjecture.

\textbf{Conjecture:} Let
$n=2pq$. Then the degree of  lowest degree polynomial in $H_n$ with at least three terms is
$min\{n-2q, n-2p-q+1\}$.\\
Hence  lowest degree polynomial(s) in $H_n$ with at least three terms is (are)

$$ \begin{cases} g_1(x), & {\mbox{ if }} \; 2p<q+1 \\
 \Phi_n^T(x), &  {\mbox{ if }} \; 2p>q+1\\
\mbox{both}\; \Phi_n^T(x)\;\mbox{and}\;g_1(x), & {\mbox{ if }} \; 2p=q+1.
\end{cases}$$

For example,  if $n=30$, then  $2p=q+1$ and
$\Phi_n^T(x)=x^{20}+x^{19}+x^{18}+x^8+x^7+x+1$
whereas $g_1(x)=x^{20}+x^{10}+1$.

\begin{lemma}\label{lem:s}
Let $n$ be an even positive integer. Let  $f(x)\in H_n$, where $|E_f|\geq 3$  and $\deg(f(x))<\frac{n}{2}+\varphi(n)$.
 Then there exists $s\in E_f$ such that $\varphi(n)\leq s< \frac{n}{2}$.
\end{lemma}

\begin{proof}
From the definition of $H_n$, it is clear that $\frac{n}{2}\notin E_f$. So we can write $f(x)=f_1(x)+f_2(x)$
in such a way that $min(E_{f_1})>\frac{n}{2}$ and
$max(E_{f_2})<\frac{n}{2}$. Then
$\Phi_n(x)|f_2(x)-x^{-(n/2)}f_1(x)$. From the given hypothesis
$\deg(x^{-(n/2)}f_1(x))<\varphi(n)$, consequently
 $max(E_{f_2})\geq \varphi(n)$. Hence the result follows.
\end{proof}

\end{document}